%% file: surveyLSymarXiv.tex
\documentclass{amsart}
\usepackage{tabmac}
\usepackage{epsfig}
\usepackage{pstricks}
\usepackage{pst-plot}
\usepackage{pdflscape}
\usepackage{color}
\usepackage{epsfig,amsthm,amsopn,amsfonts,amssymb,tabmac,pb-diagram}
\usepackage[hypertex]{hyperref}
\usepackage{verbatim}

\newtheorem{theorem}{Theorem}[section]
\newtheorem{thm}[theorem]{Theorem}
\newtheorem{lemma}[theorem]{Lemma}

\newtheorem{proposition}[theorem]{Proposition}
\newtheorem{prop}[theorem]{Proposition}

\newtheorem{problem}[theorem]{Problem}

\theoremstyle{definition}
\newtheorem{remark}[theorem]{Remark}
\newtheorem{example}[theorem]{Example}

\def\x{\bar x}
\def\R{\mathbb{R}}
\def\pp{\tilde p}
\def\C{\mathbb{C}}
\def\D{\bar D}
\def\Z{\mathbb{Z}}
\def\Q{\mathbb{Q}}
\def\Sym{\mathrm{Sym}}
\def\LSym{\mathrm{LSym}}
\def\MSym{\mathrm{MSym}}
\def\prim{\mathrm{prim}}
\def\trop{\mathrm{trop}}
\def\wt{\mathrm{wt}}
\def\het{\mathrm{ht}}

\newcommand\defn[1]{{\bf #1}}

\newcommand{\hako}{ 
\thicklines
\put(0,0){\line(0,1){0.8}}
\put(0,0){\line(1,0){0.8}}
\put(0.8,0){\line(0,1){0.8}}
\thinlines
}

\newcommand{\tama}{ 
\color{black}
\put(0.14,0.4){\circle*{0.6}}
\color{black}
\put(0.14,0.4){\circle{0.6}}
}

\begin{document}

\author{Thomas Lam}
\title{Loop symmetric functions and factorizing matrix polynomials}
\email{tfylam@umich.edu}
\address{Department of Mathematics,
University of Michigan, 530 Church St., Ann Arbor, MI 48109 USA}
\thanks{T.L. is supported by NSF grant DMS-0901111, and by a Sloan Fellowship.}

\begin{abstract}
These are notes for my talk at ICCM 2010, Beijing.  We survey some results, obtained jointly with Pavlo Pylyavskyy, concerning the ring of loop symmetric functions.  Motivations from networks on surfaces, total positivity, crystal graphs, and discrete integrable systems are discussed.
\end{abstract}
\maketitle
\section{Introduction}
The ring of symmetric functions naturally occurs when studying the problem of factorizing a polynomial into linear factors.  If one instead considers factorizations of matrix polynomials, one obtains the ring $\LSym$ of {\it loop symmetric functions}.  Whereas the ring of symmetric functions is the ring of invariants of a natural symmetric group action on the polynomial ring, the ring of loop symmetric functions is the ring of polynomial invariants of a birational symmetric group action.

The first half of this article (Sections \ref{sec:sym}-\ref{sec:loopschur}) is a brief introduction to $\LSym$, and certain distinguished elements called {\it loop Schur functions}.

Our work with Pavlo Pylyavskyy on loop symmetric functions began with our study of the theory of total positivity of loop groups \cite{TP1,LPsym}.  This connection is discussed in Section \ref{sec:TP}, and the character theory of the infinite symmetric group is given as motivation.

The birational $S_m$-action for which $\LSym$ is the invariants arises as the birational $R$-matrix of certain affine geometric crystals \cite{BK,E,KNO}.  We explain this connection at the level of combinatorial crystals in Section \ref{sec:crystal}, and show how loop Schur functions can be applied to the study of the {\it energy function} of affine crystals \cite{LPenergy}.   The energy function (for the affine crystals in this article) is equivalent to Lascoux-Sch\"{u}tzenberger's {\it cocharge statistic} on tableaux, and it is the connection to the latter that we describe.

This birational symmetric group action also leads to discrete integrable systems.  We describe Takahashi-Satsuma's box-ball system \cite{TS} in Section \ref{sec:bb}, and explain the relation, due to Hatayama, Hikami, Inoue, Kuniba, Takagi and Tokihiro \cite{HHIKTT}.  In this context $\LSym$ can be viewed as integrals of motions of these dynamical systems.

In Section \ref{sec:net}, we explain the connection between $\LSym$ and networks on the cylinder \cite{LPnet}.  On the one hand, this theory is a powerful technique  for establishing properties of the birational symmetric group action.  On the other hand, these topological networks lead to a broad collection of generalizations.

The birational $S_m$-action of this paper also arises in work of Noumi and Yamada \cite{NY} on discrete Painlev\'{e} systems, though we shall not discuss this connection.

\medskip

My results reported upon in this article are all joint work with Pavlo Pylyavskyy, and I am grateful to him for our long collaboration.

\section{Symmetric polynomials}\label{sec:sym}

\subsection{Via factorizing polynomials}
Let $p(t) = (1+x_1t)(1+x_2t) \cdots (1+x_mt)$ be a polynomial with roots $-x_1^{-1},-x_2^{-1},\ldots,-x_m^{-1}$.  The \defn{elementary symmetric polynomials} $e_1(x_1,x_2,\ldots,x_m)$,  $e_2(x_1,x_2,\ldots,x_m)$, $\ldots$, $e_m(x_1,x_2,\ldots,x_m)$ are defined by the equation
$$
p(t) = 1+ e_1\, t^{1} + e_2 \, t^{2} + \cdots + e_m\,t^m.
$$

Thus
\begin{align*}
e_1(x_1,\ldots,x_m) &= x_1 + x_2 + \cdots + x_m \\
e_2(x_1,\ldots,x_m) &= x_1x_2 + x_1x_3 + \cdots +x_1x_m+x_2x_3 + \cdots +x_{m-1}x_m
\end{align*}
and so on.

\subsection{Via an action of the symmetric group}
Let the symmetric group $S_m$ act via algebra automorphisms on $\Z[x_1,x_2,\ldots,x_m]$ by the formula
\begin{equation}\label{eq:sym}
w \cdot f(x_1,x_2,\ldots,x_m) = f(x_{w(1)},x_{w(2)},\ldots,x_{w(m)}).
\end{equation}

The following ``fundamental theorem of symmetric polynomials'' is often attributed to Newton.
\begin{theorem}\label{T:fund}
The elementary symmetric polynomials $e_1,e_2,\ldots,e_m$ are algebraically independent generators of the ring of invariants $\Z[x_1,x_2,\ldots,x_m]^{S_m}$.
\end{theorem}

We shall denote the ring of Theorem \ref{T:fund} by $\Sym_m$, and call it the \defn{ring of symmetric polynomials in $m$ variables}.  The inverse limit of the $\Sym_m$, is the \defn{ring of symmetric functions} in infinitely many variables, denoted $\Sym$.

\subsection{Via Galois Theory}
The symmetric group $S_m$ also acts on the field $\Q(x_1,x_2,\ldots,x_m)$ of rational functions in $m$ indeterminates via the same formula \eqref{eq:sym}.  In this context one may think of $S_m$ as the Galois group of the Galois extension $\Q(x_1,x_2,\ldots,x_m)/\Q(e_1,e_2,\ldots,e_m)$, which has degree $m!$.

This point of view has the weakness that one has to work with fields, rather than rings.

\section{Schur polynomials}
We refer the reader to \cite{EC2} for more concerning the material of this section.
The ring $\Sym_m$ has a $\Z$-basis consisting of the Schur polynomials $s_\lambda(x_1,x_2,\ldots,x_m)$, where $\lambda = \lambda_1 \geq \lambda_2 \geq \cdots \geq \lambda_m \geq 0$ is a partition with at most $m$ parts.

\subsection{As a generating function of tableaux}
If $\mu \subset \lambda$ as Young diagrams, then we can define the skew Schur polynomial
$$
s_{\lambda/\mu}(x_1,x_2,\ldots,x_m) = \sum_T x^{\wt(T)}
$$
where the sum is over all skew semistandard Young tableaux (SSYT) with shape $\lambda/\mu$, filled with the letters $1,2,\ldots,m$ and $\wt(T)$ denotes the weight of the tableau $T$.  For example,
$$
T = \tableau[sY]{\bl &\bl &\bl&1&1&3&5&5\\ \bl&1&2&2&2&5&6 \\4&4&6}
$$
is a SSYT with shape $(8,7,3)/(3,1)$ and $x^{\wt(T)} = x_1^3 x_2^3 x_3 x_4^2 x_5^3 x_6^2$.  When $\mu = (0)$ then $s_{\lambda/(0)} = s_\lambda$ is a Schur polynomial.

\subsection{As a Jacobi-Trudi determinant}
If $\lambda$ is a partition, then $\lambda'=(\lambda'_1,\ldots,\lambda'_\ell)$ denotes the conjugate partition, where $\lambda'_i$ is the number of boxes in the $i$-th column of $\lambda$.

\begin{theorem}[Jacobi-Trudi formula]\label{T:JT}
The skew Schur polynomial $s_{\lambda/\mu}$ has the following determinantal expression in terms of elementary symmetric polynomials:
$$
s_{\lambda/\mu} = \det(e_{\lambda'_i-i+j-\mu'_j})_{i,j=1}^m.
$$
\end{theorem}

In particular, the skew Schur polynomials $s_{\lambda/\mu}(x_1,\ldots,x_m)$ lie in $\Sym_m$.

\subsection{As a ratio of alternants}
Let $\alpha=(\alpha_1,\alpha_2,\ldots,\alpha_m)$ be a sequence of nonnegative integers.  We define the \defn{alternant}
$$
a_\alpha = \sum_{w \in S_m} (-1)^w \, w(x^\alpha)
$$
where $(-1)^w$ denotes the sign of $w$.  The alternant can also be expressed as the determinant $a_\alpha = \det(x_i^{\alpha_j})_{i,j=1}^m$.  Define the staircase partition $\delta_{m} := (m,m-1,\ldots,0)$.

\begin{theorem}\label{T:alternant}
For $\lambda$ with less than or equal to $m$ parts,
$$
s_\lambda(x_1,\ldots,x_m) = \frac{a_{\lambda+\delta_{m-1}}}{a_{\delta_{m-1}}}.
$$
\end{theorem}


\section{Loop symmetric polynomials}

\label{sec:lsym}
Fix an integer $n \geq 1$, which will usually be suppressed in our notation.  We shall define a ring $\LSym_m$ of \defn{loop symmetric polynomials in $m$ sets of variables}.  When $n = 1$, the ring $\LSym_m$ reduces to $\Sym_m$.

\subsection{Via solving matrix equations}
Suppose one is given a $n\times n$ matrix $P(t)$ with coefficients which are polynomials in $t$.  Is it possible to factorize $P(t)$ into a product of ``linear'' matrix factors?  In the following we will formulate precisely one version of this problem.

Define the $n \times n$ \defn{whirl matrix}
$$
M(a_1,a_2,\ldots,a_n) = \left(\begin{array}{cccccc}1 & a_1 & 0 & \cdots & 0 &0  \\
0&1&a_2&\cdots&0&0 \\
0&0&1&\cdots&0&0 \\
\vdots& \vdots & \vdots &\vdots &\vdots &\vdots \\
0&0&0&\cdots&1&a_{n-1} \\
a_nt &0 &0 &\cdots &0&1  \end{array} \right)
$$
For $n = 1$, we define $M(a) = 1+at$.  We think of $M(a_1,a_2,\ldots,a_n)$ as a polynomial (in $t$)-valued matrix, which in particular is linear.  Every matrix entry of $M$ has degree less than or equal to $1$.

Now let $P(t)$ be a $n \times n$ matrix with polynomial coefficients.  Can one factorize
\begin{equation}\label{E:whirlfact}
P(t) = M(x_1^{(1)},\ldots,x_1^{(n)}) \cdots M(x_m^{(1)},\ldots,x_m^{(n)})
\end{equation}
for suitable ``roots'' $x_i^{(k)}$?  There is a natural condition to enforce: each matrix $M(x_1^{(1)},\ldots,x_n^{(1)})$ becomes an upper triangular matrix with $1$'s on the diagonal when $t = 0$.  So to obtain a solution, we must at least assume that $P(0)$ is an upper triangular matrix with $1$'s on the diagonal (uni-upper triangular).  This is analogous to the fact that the constant term of $p(t) = (1+x_1t)(1+x_2t)\cdots(1+x_mt)$ must be $1$.

\begin{theorem}[{\cite[Theorem 4.1]{LPcrystal}}]
Suppose $P(t)$ is a $n \times n$ matrix with entries in $\C[t]$, such that $P(0)$ is uni-upper triangular.  If $P(t)$ is generic, then it has a factorization of the form \eqref{E:whirlfact}.
\end{theorem}

Given a factorized matrix polynomial $P(t) = M(x_1^{(1)},\ldots,x_1^{(n)}) \cdots M(x_m^{(1)},\ldots,x_m^{(n)})$, we define the \defn{loop elementary symmetric functions} $\{e_i^{(k)} \mid i =1,\ldots,n \;\text{and}\;k=1,\ldots,m\}$ by the matrix equation \

\noindent
\resizebox{\hsize}{!}{
$
P(t) =  \left(\begin{array}{lllll}1+e_n^{(1)}t+ \cdots & e_1^{(1)}+e_{n+1}^{(1)}t + \cdots  & \cdots & e_{n-2}^{(1)}+e_{2n-2}^{(1)}t + \cdots & e_{n-1}^{(1)}+e_{2n-1}^{(1)}t + \cdots \\
e_{n-1}^{(2)}t+e_{2n-1}^{(2)}t^2+ \cdots &1+e_n^{(2)}t+ \cdots  & \cdots& e_{n-3}^{(2)}+e_{2n-3}^{(2)}t + \cdots & e_{n-2}^{(2)}+e_{2n-2}^{(2)}t + \cdots\\
\vdots& \vdots & \vdots &\vdots &\vdots \\
e_2^{(n-1)}t+e_{n+2}^{(n-1)}t^2 + \cdots&e_3^{(n-1)}t+e_{n+3}^{(n-1)}t^2 + \cdots& \cdots&1+e_n^{(n-1)}t+ \cdots&e_1^{(n-1)}+e_{n+1}^{(n-1)}t +\cdots\\
e_1^{(n)}t+e_{n+1}^{(n)}t^2+\cdots &e_2^{(n)}t+e_{n+2}^{(n)}t^2 + \cdots  &\cdots &e_{n-1}^{(n)}t + e_{2n-1}^{(n)}t^2+\cdots &1+e_n^{(n)}t+ \cdots  \end{array} \right)
$
}
\

\noindent
For example, let us take $n = 2$ and $m = 3$.  One has
\begin{align*}
&\left(\begin{array}{cc}1&x_1^{(1)}\\x_1^{(2)}t&1 \end{array}\right)\left(\begin{array}{cc}1&x_2^{(1)}\\x_2^{(2)}t&1 \end{array}\right)
\left(\begin{array}{cc}1&x_3^{(1)}\\x_3^{(2)}t&1 \end{array}\right)
\\
&=
\left(\begin{array}{cc}1+(x_1^{(1)}x_2^{(2)}+x_2^{(1)}x_3^{(2)}+x_1^{(1)}x_3^{(2)})t&
x_1^{(1)}+x_2^{(1)}+x_3^{(1)} + x_1^{(1)}x_2^{(2)}x_3^{(1)}t \\
(x_1^{(2)}+x_2^{(2)}+x_3^{(2)})t + x_1^{(2)}x_2^{(1)}x_3^{(2)}t^2 & 1+(x_1^{(2)}x_2^{(1)}+x_1^{(2)}x_3^{(1)}+x_2^{(2)}x_3^{(1)})t\end{array}\right)
\end{align*}
so that
$$
\begin{array}{ll}
e_1^{(1)} = x_1^{(1)}+x_2^{(1)}+x_3^{(1)}  &e_1^{(2)} = x_1^{(2)}+x_2^{(2)}+x_3^{(2)}\\
e_2^{(1)} =x_1^{(1)}x_2^{(2)}+x_2^{(1)}x_3^{(2)}+x_1^{(1)}x_3^{(2)} & e_2^{(2)}=x_1^{(2)}x_2^{(1)}+x_1^{(2)}x_3^{(1)}+x_2^{(2)}x_3^{(1)}\\
e_3^{(1)}=x_1^{(1)}x_2^{(2)}x_3^{(1)} & e_3^{(2)}=x_1^{(2)}x_2^{(1)}x_3^{(2)}
\end{array}
$$
Note that $e_k^{(r)}$ is a homogeneous polynomial of degree $i$.  An explicit formula for $e_k^{(r)}$ is
\begin{align*}
e_k^{(r)}(x_1,x_2,\ldots,x_m) &= \sum_{1 \leq i_1 < i_2 < \cdots < i_k\leq m} x_{i_1}^{(r)} x_{i_2}^{(r+1)} \cdots x_{i_k}^{(r+k-1)}
\end{align*}
By convention, $e_k^{(r)} = 0$ for $k < 0$, and $e_0^{(r)}=1$.  We call the upper index the \defn{color}.  When all $n$ colors are identified, that is $x_i^{(s)}=x_i^{(s')}$ for all $i$ and $s,s'\in \Z/n\Z$, these polynomials specialize to the usual elementary symmetric polynomials.  We denote the subring of  $\Z[x_1^{(1)},x_1^{(2)},\ldots,x_1^{(n)},x_2^{(1)},\ldots,x_m^{(n)}]$ generated by the loop elementary symmetric polynomials $e_i^{(k)}$ by $\LSym_m$, the ring of \defn{loop symmetric polynomials in $m$ sets of variables}.  Our loop elementary symmetric polynomials (with some change in labeling of variables) appeared previously in work of Yamada \cite{Y}.

\begin{remark}
In \cite{TP1} we distinguished two different rings of loop symmetric functions, one called the {\it curl ring}, and the other the {\it whirl ring}.  The two are related by negation of the upper index $x_i^{(r)} \mapsto x_i^{(-r)}$, where as usual the upper index is taken modulo $n$.
\end{remark}

\subsection{Via the birational $R$-matrix}\label{sec:R}
The whirl matrices $M(x)$ do not in general commute with each other for $n > 1$, but they do satisfy a commutation relation.

For two $n$-tuples $(x^{(1)},\ldots,x^{(n)})$ and $(y^{(1)},\ldots,y^{(n)})$ of variables, we denote
$$
\kappa_i(x,y) = \sum_{j=i}^{i+n-1} \prod_{k=i+1}^j y^{(k)}\; \prod_{k=j+1}^{i+n-1} x^{(k)}
$$
where here and elsewhere, upper indices are always taken modulo $n$.  For example, for $n=3$, one has
$$
\kappa_1(x,y) = x^{(2)}x^{(3)}+y^{(2)}x^{(3)}+y^{(2)}y^{(3)}.
$$

Now define a rational map 
$$s: \Q(x^{(1)},\ldots,x^{(n)},y^{(1)},\ldots,y^{(n)}) \to \Q(x^{(1)},\ldots,x^{(n)},y^{(1)},\ldots,y^{(n)})$$ 
by
\begin{equation}\label{E:s}
s(x^{(i)})=y^{(i+1)}\frac{\kappa_{i+1}(x,y)}{\kappa_i(x,y)} \ \ \ \ \text{and} \ \ \ \  s(y^{(i)})=x^{(i-1)}\frac{\kappa_{i-1}(x,y)}{\kappa_i(x,y)}.
\end{equation}

If we are given $m$ sets of variables $\{(x_i^{(1)},\ldots,x_i^{(n)}) \mid 1 \leq i \leq m\}$, then we let $s_k$ denote the rational map acting on $\Q(x^{(i)}_j)$, fixing the sets of variables $x_i$ for $i \neq k, k+1$, and then applying $s$ to the variables $x_k, x_{k+1}$.

\begin{theorem}\label{T:R} \
\begin{enumerate}
\item
The rational maps $s_j$ generate a birational action of $S_m$ on $\Q(x^{(i)}_j)$.
\item
Let $w \in S_m$ act on $\Q(x^{(i)}_j)$ as in (1).  Then
$$
M(x_1)M(x_2)\cdots M(x_m) = M(w(x_1))M(w(x_2))\cdots M(w(x_m)).
$$
\end{enumerate}
\end{theorem}

Note that if $n = 1$, then this birational action of $S_m$ is just the usual action of $S_m$ on $\Q(x_1,\ldots,x_m)$ of Section \ref{sec:sym}.  The birational action of $S_m$, and the proof of Theorem \ref{T:R}(1) was established independently in a number of different contexts \cite{E,Ki,NY,TP1,LPnet}.  Theorem \ref{T:R}(2) is phrased in terms of ``M-matrices'' in the context of affine geometric crystals \cite{KNO}.

The following result will be established in \cite{LPsym}.  Its proof is significantly harder than the fundamental theorem of symmetric functions.
\begin{theorem}[Fundamental theorem of loop symmetric functions]
\label{T:fundlsym}
The ring of loop symmetric polynomials is exactly the polynomial invariants of the birational $S_m$-action:
$$\LSym_m = \Q(x^{(i)}_j)^{S_m} \cap \Z[x^{(i)}_j].$$
Furthermore, the loop elementary symmetric polynomials are algebraically independent generators of $\LSym_m$.
\end{theorem}

We note that the definition of $\LSym_m$ is compatible with restriction of variables: setting $x_m = 0$ gives $\LSym_{m-1}$ from $\LSym_m$.  Thus it is possible to define a ring $\LSym$, in infinitely many variables, as the inverse limit of $\LSym_m$.  We call $\LSym$ the ring of \defn{loop symmetric functions}.

\begin{example}
Let $n = 3$ and $m = 2$, and set $x = x_1$ and $y = x_2$.  Then
\begin{align*}
s(e^{(1)}_1) &= s(x^{(1)}+y^{(1)}) \\
&= y^{(2)}\frac{\kappa_2}{\kappa_1} + x^{(3)}\frac{\kappa_3}{\kappa_1} \\
&= \frac{y^{(2)}(x^{(3)}x^{(1)}+y^{(3)}x^{(1)}+y^{(3)}y^{(1)}) + x^{(3)}(x^{(1)}x^{(2)}+y^{(1)}x^{(2)}+y^{(1)}y^{(2)})}{x^{(2)}x^{(3)}+y^{(2)}x^{(3)}+y^{(2)}y^{(3)}}\\
&= x^{(1)}_1+y^{(1)}_1 = e^{(1)}_1.
\end{align*}
\end{example}

\subsection{Galois groups of matrix polynomials}
It can be deduced from the preceding discussion that $\Q(x_i^{(j)})/\Q(x_i^{(j)})^{S_m}$ is a Galois extension of degree $m!$, with Galois group $S_m$.  (Indeed, this is significantly easier than Theorem \ref{T:fundlsym}, see \cite[Theorem 4.1]{LPcrystal}.)  Thus one can think of $S_m$ as the Galois group of the matrix polynomial $P(t)$, where the coefficients of the polynomial are the ``variables'' $e^{(i)}_j$.  Perhaps it would be interesting to ask: what can be said about the ``Galois groups'' of matrix polynomials with coefficients taking values in $\Z$ or $\Q$?

\section{Loop skew Schur functions}
\label{sec:loopschur}
\subsection{As a generating function of tableaux}
A square $s = (i,j)$ in the $i$-th row and $j$-th column has \defn{content} $c(s)=i-j$.  We caution that our notion of content is the negative of the usual one.  Let $\rho/\nu$ be a skew shape.   Recall that a semistandard Young tableaux $T$ with shape
$\rho/\nu$ is a filling of each square $s \in \rho/\nu$ with an
integer $T(s) \in \Z_{> 0}$ so that the rows are weakly-increasing, and columns are increasing.  For $r\in \Z/n\Z$, the \defn{$r$-weight} $x^{\wt^{(r)}(T)}$ of a tableaux $T$ is given by $x^{\wt^{(r)}(T)} = \prod_{s \in \rho/\nu}x_{T(s)}^{(c(s)+r)}$.

We shall draw our shapes and tableaux in English notation:
$$
\tableau[sY]{\bl \circ&\bl \circ&\bl \circ&\bl \circ&\bl \circ&&&&\\\bl \circ&\bl \circ&\bl \circ&&&&&\\\bl \circ&\bl \circ&\bl \circ&&&} \qquad
\tableau[sY]{\bl &\bl&1&1&1&3\\1&2&2&3&4\\3&3&4}
$$
%

\noindent For $n = 3$ the $0$-weight of the above tableau is $(x_1^{(1)})^2 (x_3^{(1)})^3
x_1^{(2)} x_2^{(2)} x_3^{(2)} x_1^{(3)} x_2^{(3)} (x_4^{(3)})^2.$
We define the \defn{loop (skew) Schur function} by
$$
s^{(r)}_{\lambda/\mu}({x}) = \sum_{T} x^{\wt^{(r)}(T)}
$$
where the summation is over all semistandard Young tableaux of
(skew) shape $\lambda/\mu$.

\begin{example}
Let $n = 2$.  Then
\begin{align*}
s^{(1)}_{2,1}(x_1,x_2,x_3) = & x_1^{(1)}x_1^{(2)}x_2^{(2)} + x_1^{(1)}x_2^{(2)}x_2^{(2)} + x_1^{(1)}x_2^{(2)}x_3^{(2)} +
x_1^{(1)}x_3^{(2)}x_2^{(2)} + \\
&x_1^{(1)}x_1^{(2)}x_3^{(2)} +
x_2^{(1)}x_2^{(2)}x_3^{(2)} +
x_1^{(1)}x_3^{(2)}x_3^{(2)} +
x_2^{(1)}x_3^{(2)}x_3^{(2)}
\end{align*}
corresponding to the tableaux
$$
\tableau[sY]{1&1\\2} \qquad \tableau[sY]{1&2\\2} \qquad \tableau[sY]{1&2\\3} \qquad \tableau[sY]{1&3\\2}
$$
$$
\tableau[sY]{1&1\\3} \qquad \tableau[sY]{2&2\\3} \qquad \tableau[sY]{1&3\\3} \qquad \tableau[sY]{2&3\\3}
$$
where the upper left corner has content 0, and so gives a $1$-weight with color $1$.  Setting $x_i^{(1)}=x_i^{(2)} = x_i$ gives the usual Schur polynomial $s_{2,1}(x_1,x_2,x_3)$.
\end{example}

Perhaps surprisingly, the loop Schur functions do not form a basis for $\LSym$.  By Theorem \ref{T:fundlsym}, monomials in the loop elementary symmetric functions will form a basis for $\LSym$, but finding the analogue of monomial symmetric functions is not straightforward.   The following problems seem central.

\begin{problem}
Find a Schur-like basis for $\LSym$.
\end{problem}

\begin{problem}
Find a monomial-like basis for $\LSym$.
\end{problem}

Presumably any ``loop monomial symmetric function'' would be nonnegative in terms of monomials of the variables, but be minimally so in some sense.

\subsection{As a determinant}
We have the following analogue of the Jacobi-Trudi formula.

\begin{thm}[{\cite[Theorem 7.6]{TP1}}] \label{T:loopJT} The loop skew Schur function has the following determinantal expression in terms of loop elementary symmetric functions: $$s^{(r)}_{\lambda'/\mu'} = \det(e_{\lambda_i-\mu_j-i+j}^{(r-j+1+\mu_j)}).$$
Thus $s^{(r)}_{\lambda/\mu} \in \LSym$.
\end{thm}

The proof of Theorem \ref{T:loopJT} is essentially the same as that of Theorem \ref{T:JT}.

\begin{remark}
The Jacobi-Trudi formula puts loop Schur functions into the context of Macdonald's ninth variation of Schur functions.  The latter have a similar Jacobi-Trudi formula \cite{NNSY}, for which the upper index is not cyclic.
\end{remark}

\subsection{As a ratio of alternants}
Suppose we are given an $n \times m$ array $\{x_i^{(j)} \mid 1\leq i \leq m \;\; \text{and} \;\; j \in \Z/n\Z\}$ of variables.  Let $r \in \Z/n\Z$.  Given a decreasing sequence $\alpha_1 > \alpha_2 > \cdots > \alpha_n \geq 0$ of nonnegative integers, we define the \defn{loop alternant} $a^{(r)}_\alpha$ as the $m \times m$ determinant
$$
\det\left(t_{m-j+1,m}(x_m^{(r)}x_m^{(r-1)} \cdots x_m^{(r-\alpha_i+1)})\right)_{i,j=1}^m
$$
where $t_{a,b}$ denotes the transposition exchanging $a$ and $b$, acting via the birational action of $S_m$.  The following result \cite{LPsym} generalizes Theorem \ref{T:alternant}\footnote{In an earlier version of this paper, this result was stated as a conjecture.  We thank Greg Anderson for the idea of the proof.}.

\begin{thm}\label{thm:alt}
For $\lambda$ with less than or equal to $m$ parts, the loop Schur function $s_\lambda^{(r-1)}$ can be expressed as
$$
s_\lambda^{(r-1)}(x_1,\ldots,x_m) = a^{(r)}_{\lambda+\delta}/a^{(r)}_{\delta}.
$$
\end{thm}

\begin{example}
Let $\lambda = (2,1)$, $n = 3$ and $m = 2$, denoting the variables by $x= x_1$ and $y = x_2$.  Then
\begin{align*}
a^{(1)}_{3,1} &= \det\left(\begin{array}{cc} y^{(1)}y^{(3)}y^{(2)} &x^{(1)}x^{(3)}x^{(2)}\\
y^{(1)} & x^{(3)}\kappa_3/\kappa_1  \end{array}\right)\\
a^{(1)}_{1,0} &= \det\left(\begin{array}{cc} y^{(1)}&x^{(3)}\kappa_3/\kappa_1  \\ 1&1
\end{array}\right)
\end{align*}
So
{\small
\begin{align*}
&a^{(1)}_{3,1}/a^{(1)}_{1,0}\\
&= \frac{y^{(1)}y^{(2)}y^{(3)}x^{(3)}(x^{(1)}x^{(2)}+y^{(1)}x^{(2)}+y^{(1)}y^{(2)}) - x^{(1)}x^{(2)}x^{(3)}y^{(1)}(x^{(2)}x^{(3)} + y^{(2)}x^{(3)}+y^{(2)}y^{(3)})}{x^{(1)}x^{(2)}x^{(3)}-y^{(1)}y^{(2)}y^{(3)}}\\
&=x^{(3)}y^{(1)}x^{(2)}+x^{(3)}y^{(1)}y^{(2)}\\
&= s^{(3)}_{2,1}(x,y).
\end{align*}
}
\end{example}

\subsection{$\LSym$ as a Hopf algebra}
The ring of symmetric functions has a Hopf algebra structure, where comultiplication is given by $\Delta(e_i) = \sum_{j=0}^i e_j \otimes e_{i-j}$, and the antipode is given by $S(e_i) = (-1)^{i} h_i$, where $h_i$ is the complete homogeneous symmetric function.  We now briefly describe the Hopf structure of $\LSym$.

The comultiplication of $\LSym$ is given by
$$
\Delta(e_i^{(k)}) = \sum_{j=0}^i e_j^{(k)} \otimes e_{i-j}^{(k+j)}
$$
where by convention $e_0^{(k)}=1$ for every $k$.  Note that $\LSym$ is not co-commutative.  The antipode is given by
$$
S(e_i^{(k)}) = s_{(i)}^{(k+i-1)}.
$$
Clearly $\LSym$ is graded, with $\deg(e_i^{(k)}) = i$.

\begin{prop}
The above structures define a graded Hopf algebra $\LSym$.
\end{prop}

Let $\MSym$, the \defn{Hopf algebra of matrix symmetric functions} be the graded dual Hopf algebra to $\LSym$.  By general theory, one knows that $\MSym_\Q = \MSym \otimes_\Z \Q$
is generated (as an algebra) by its primitive subspace $\prim$.  In fact $\MSym_\Q$ is the universal enveloping algebra of the Lie algebra $\prim$.

\begin{prop}
The primitive subspace $\prim$ of $\MSym_\Q$ is the the Lie-algebra of
infinite $\Z \times \Z$ indexed matrices, which are $n$-periodic,
strictly upper-triangular, and which contains finitely many non-zero
diagonals.  The grading on $\prim$ is such that the $i$-th diagonal
above the main diagonal has degree $i$.
\end{prop}

\subsection{Loop powersum symmetric functions, and a loop Murnaghan-Nakayama rule}
For $k \geq 1$, define the \defn{loop powersum symmetric polynomials}
$$
\pp_k(x_1,x_2,\ldots,x_m) = (x_1^{(1)}x_1^{(2)} \cdots x_1^{(n)})^k + (x_2^{(1)}x_2^{(2)} \cdots x_2^{(n)})^k + \cdots
(x_m^{(1)}x_m^{(2)} \cdots x_m^{(n)})^k
$$
so that $\pp_k$ is homogeneous of degree $kn$.  The $\pp_k(x_1,x_2,\ldots)$ are primitive elements of $\LSym$, though even the fact that they lie in $\LSym$ is not completely obvious.  Recall that a ribbon is a connected skew shape not containing any $2\times 2$ square.

\begin{theorem}[{\cite{LPsym}}]
We have
$$
\pp_k s_\lambda^{(r)} = \sum_{\mu} (-1)^{\het(\mu/\lambda)} \, s_\mu^{(r)}
$$
where the summation is over all ribbons $\mu/\lambda$ of size $kn$, and $\het(\mu/\lambda)$ denotes the number of rows in $\mu/\lambda$, minus 1.
\end{theorem}

Thus the usual Murnaghan-Nakayama rule \cite{EC2} holds for the loop powersum and loop Schur functions.

\section{Networks on cylinders and the birational $S_m$-action}\label{sec:net}
\subsection{Local description of  $S_m$-action}
The fact that the transformation \eqref{E:s} satisfies the braid relation (Theorem \ref{T:R}(1)) is not obvious.  Here we sketch an approach based on a topological model of networks developed in \cite{LPnet}.

\begin{figure}[h!]
    \begin{center}
    \input{wire181shift.pstex_t}
    \end{center}
    \caption{The birational action $(x,y) \mapsto (x',y')$.}
    \label{fig:wire181}
\end{figure}
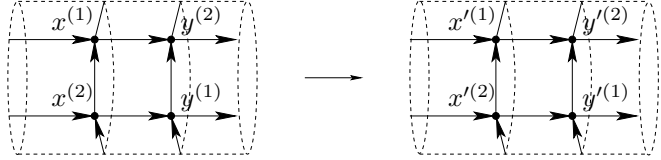

Let us first arrange our variables $x_i^{(j)}$ on a network embedded into a cylinder.  The variables are located on the vertices of this network, that consists of $n$ horizontal wires connecting the boundaries of the cylinder, and $m$ loops which go around the cylinder.  The case $n = m =2$ is illustrated in Figure \ref{fig:wire181}.  The birational action $s: \Q(x_i^{(j)}) \to \Q(x_i^{(j)})$ will be expressed in terms of certain {\it local transformations} of the network.

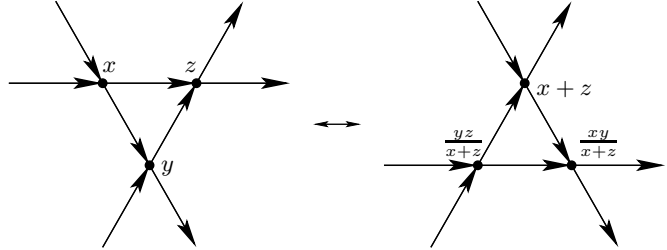
\begin{figure}[h!]
    \begin{center}
    \input{wire8shrink.pstex_t}
    \end{center}
    \caption{Yang-Baxter move with transformation of vertex weights shown.}
    \label{fig:wire8}
\end{figure}

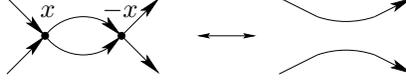
\begin{figure}[h!]
    \begin{center}
    \input{wire11a.pstex_t}
    \end{center}
    \caption{Adding or removing a crossing.}
    \label{fig:wire11}
\end{figure}

We need the following two local transformations (more are discussed in \cite{LPnet}), the \defn{Yang-Baxter relation} (Figure \ref{fig:wire8}) and the \defn{crossing creation/removal rule} (Figure \ref{fig:wire11}).  The parameters in this Yang-Baxter rule naturally occur in Lusztig's study of total positivity in reductive groups \cite{Lus}.

To calculate the birational action $s$, we first use the crossing addition rule to create a pair of crossings with weights $p$ and $-p$.  Using the Yang-Baxter move, we may push one of these crossings through to the horizontal wires, until it loops around the cylinder.  In effect, this has swapped the two vertical loops.  For a unique choice of the weight $p$, the weight that comes out on the other side after passing through all the horizontal wires is also $p$; and so one can use the crossing removal rule the pair of crossings with weights $p$ and $-p$.

\begin{figure}[h!]
    \begin{center}
    \scalebox{0.7}{\input{wire20shift.pstex_t}}
    \end{center}
    \caption{Swapping two vertical loops by creating a pair of crossings, pushing one of them around the cylinder, and removing these crossings.}
    \label{fig:wire20}
\end{figure}
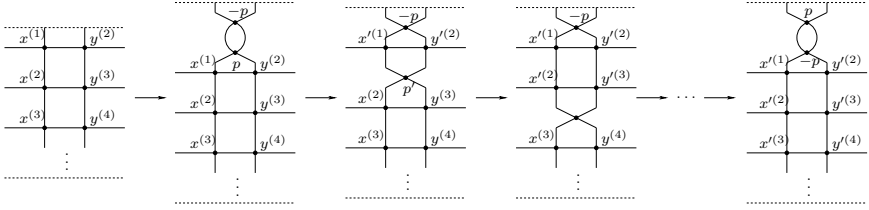

\begin{theorem}[{\cite{LPnet}}]\label{T:local}
The above local procedure calculates the birational $S_m$-action of Section \ref{sec:R}.
\end{theorem}

Using Theorem \ref{T:local}, one can prove many properties of the birational $S_m$-action, including Theorem \ref{T:R}(1), using ``local'' computations.

\subsection{Boundary and cycle measurements}

Take a network $N_{m}$ such as the one in Figure \ref{fig:wire181}, but with $m$ vertical loops and $n$ horizontal wires (as usual $n$ is suppressed from the notation).  Then the loop elementary symmetric functions can be obtained by taking certain generating functions of paths in this network.

\begin{figure}[h!]
    \begin{center}
    \input{wire2shrink.pstex_t}
    \end{center}
    \caption{The three ways a highway path can go through a vertex with weight $x$; a fragment of a highway path, contributing the factor $x z$ to the weight of the whole path.}
    \label{fig:wire2}
\end{figure}
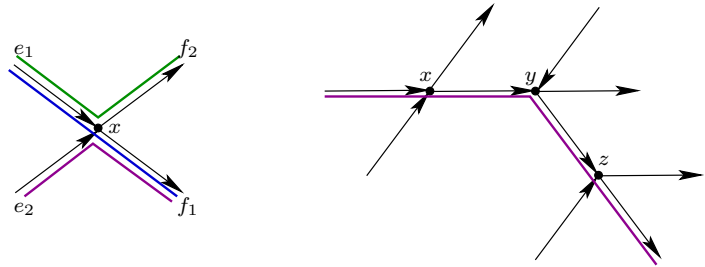

Every intersection of these networks consists of two wires intersecting transversally.  Suppose $p$ is a path.  Then at each intersection it can either go straight, or turn.  A \defn{highway path} is a path which always turns when the crossing wire crosses your current wire from the left.  The \defn{weight} of a highway path is the product of the vertex weights over all vertices where the path goes straight.  See Figure \ref{fig:wire2}.

\begin{lemma}
Loop elementary symmetric functions are exactly weight generating functions of highway paths in $N_m$ with fixed starting and ending boundary vertices, and fixed homology class.
\end{lemma}

The generating functions of the above lemma are called \defn{boundary measurements}.

\begin{example}\label{ex:loop}
Let us take the network in Figure \ref{fig:wire181}.  Then there is only one highway path which starts at the upper left vertex and ends at the upper right vertex, with weight $e_2^{(1)}(x,y) = x^{(1)}y^{(2)}$.  This is the highway path which just goes straight.  There are two highway paths from the upper left vertex to the lower right vertex, with weights $x^{(1)}$ and $y^{(1)}$ respectively, summing to $e_1^{(1)}$.  Both of these highway paths loop around the back of the cylinder.
\end{example}

In Example \ref{ex:loop}, we did not have to fix the homology class of the path.  But we would have to do so for larger $m$.

One can also define \defn{cycle measurements}.  In the case of the networks $N_m$, we should consider \defn{underway cycles}, which are cycles in $N_m$ where the condition at each intersection has left and right swapped, compared to highway paths.  Underway cycle measurements in $N_m$ give the loop powersum symmetric functions.  (There is a technical multiplicity which is involved in the definition of cycle measurements, which we omit here.)

\begin{example}
In Figure \ref{fig:wire181}, there are two underway cycles which loop around the cylinder exactly once.  These have weights $x^{(1)}x^{(2)}$ and $y^{(1)}y^{(2)}$ respectively, which sum to $\pp_1(x,y)$.  If we consider underway cycles which loop around the cylinder twice, we will get two cycles, with weights $(x^{(1)}x^{(2)})^2$ and $(y^{(1)}y^{(2)})^2$, summing to $\pp_2(x,y)$.
\end{example}

The fact that boundary and cycle measurements are preserved by the birational $S_m$-action is not a coincidence.
In \cite{LPnet}, we establish

\begin{theorem}[{\cite[Theorem 3.2]{LPnet}}]
Boundary and cycle measurements are preserved by local transformations.
\end{theorem}

Here ``local transformations'' consist of a list of moves, which includes those in Figures \ref{fig:wire8} and \ref{fig:wire11}.


\begin{remark}
The topological interpretation of $\LSym$ gives rise to many generalizations.  In \cite{LPnet}, we also prove that one obtains an action of the symmetric group when the horizontal wires do not point in the same direction (see Remark \ref{rem:BBBS} for an application of this).  However, any (valid) network on an orientable surface would give rise to a ring of invariants, where the symmetric group is replaced by a ``monodromy group'' of networks we define in \cite{LPnet}.
\end{remark}

\section{Total positivity}\label{sec:TP}
\subsection{Total positivity implies factorization}
A polynomial $p(t)$ with complex coefficients always has a factorization into a product of linear factors.  The analogous statement is not true for the matrix factorizations $$P(t) = M(x_1^{(1)},\ldots,x_1^{(n)}) \cdots M(x_m^{(1)},\ldots,x_m^{(n)})$$ we consider, but there is an elegant condition which always guarantees factorization {\it over the reals}.

Given a matrix polynomial $P(t)$, we associate to it an infinite periodic matrix $X = X_P$, as follows:
\begin{eqnarray*}
\begin{array}{c}\left(
\begin{array}{cc}
1 +4t+ 9t^2 & 2 - 2t \\
4 +6t^2  & 7+t-t^2 \end{array} \right) =  \\
\left(\begin{array}{cc}1 & 2\\ 4 & 7\end{array}\right) + t
\left(\begin{array}{cc}4 & -2\\ 0 & 1\end{array}\right) + t^2
\left(\begin{array}{cc}9 & 0 \\ 6 & -1\end{array} \right)
\end{array} &\hspace{-20pt}\rightsquigarrow& \hspace{-10pt}\left(
\begin{array}{c|cc|cc|cc|c}
\ddots & \vdots & \vdots & \vdots  & \vdots  & \vdots&\vdots \\
\hline
\dots& 4 & -2& 9& 0&0& 0&\dots \\
\dots&0&1&6&-1&0&0&\dots \\
\hline
\dots&1  & 2 & 4 & -2 & 9 & 0&\dots\\
\dots&4 & 7 &0 & 1 &6&-1 & \dots \\
\hline
\dots&0&0&1  & 2 & 4 & -2 & \dots\\
\dots&0&0& 4 & 7 &0 & 1 & \dots \\
\hline & \vdots & \vdots & \vdots  & \vdots  & \vdots & \vdots &
\ddots\end{array} \right) \\
P(t) \hspace{100pt} && \hspace{80pt} X_P
\end{eqnarray*}
In other words, $X_P$ is the block-Toeplitz matrix, where each block is equal to the coefficient of $t^i$ in $P(t)$, for some $i$.

Suppose $P(t)$ is a matrix polynomial with real coefficients.  We now define $P(t)$ to be \defn{totally nonnegative} if every minor of $X_P$ is nonnegative.

\begin{theorem}[\cite{TP1}]
Suppose $P(t)$ is a real matrix polynomial, with $P(0)$ uni-upper triangular.  Then $P(t)$ is totally nonnegative if and only if there is a factorization of $P(t)$ into a product of whirls with nonnegative real parameters.
\end{theorem}

For $n=1$, we recover a well-known criterion for a polynomial to have nonnegative real roots.

A much more difficult problem is to classify the real matrices with {\it formal power series} coefficients which are totally nonnegative.  This is studied in \cite{DMS,TP1}.  This problem can be reformulated in terms of $\LSym$.  Let $U_{\geq 0}$ denote the set of matrices $P(t)$ with real formal power series coefficients such that $P(0)$ is uni-upper triangular, and such that $P(t)$ is totally nonnegative.

\begin{proposition}\label{P:hom}
The set $U_{\geq 0}$ is in bijection with the set
$$
\{\phi:\LSym \otimes \R \to \R \mid \phi(s^{(r)}_{\lambda/\mu}) \geq 0 \}
$$
of algebra homomorphisms of $\LSym \otimes \R$ to the reals taking nonnegative values on all loop skew Schur functions.
\end{proposition}


When $n = 1$, the set $U_{\geq 0}$ is the set of \defn{totally positive functions}: formal power series $p(t) = 1+ a_1t + a_2t^2 + \cdots$ such that the infinite Toeplitz matrix
$$
\left(\begin{array}{ccccc} 1 & a_1 & a_2 & a_3 & \cdots \\
&1 & a_1 & a_2 &  \cdots \\
&&1 & a_1 &  \cdots \\
&&&1 & \cdots\\
&&&& \ddots
\end{array}\right)
$$
is totally nonnegative.

\subsection{Relation to infinite symmetric group}
The classical {\it Frobenius characteristic map} identifies the ring $\Sym$ of symmetric functions with the direct sum of the character rings of the symmetric groups $\{S_n \mid n \geq 1\}$.  It is perhaps not surprising that the set of totally positive functions has an interpretation in terms of the representation theory of the symmetric groups.

The infinite symmetric group $S_\infty$ is the union $\cup_{n \geq 1} S_n$ of all the finite symmetric groups, where $S_{n-1}$ is identified naturally with a subgroup of $S_n$.  A \defn{character} of $S_\infty$ is a complex-valued function $\chi: S_\infty \to \C$ which satisfies four conditions:
\begin{enumerate}
\item central -- $\chi$ is constant on conjugacy classes
\item positive-definite -- the matrix $(\chi(g_i^{-1} g_j))_{i,j=1}^r$ is Hermitian and nonnegative-definite for any $g_1,g_2,\ldots,g_r$
\item extremal -- $\chi$ is not a sum of two linearly independent functions satisfying (1) and (2)
\item normalized -- $\chi(1) = 1$.
\end{enumerate}
For a finite group, such a function would exactly be an irreducible character, normalized by dividing by the dimension of the character.

Let us call a totally positive function $p(t) = 1 + a_1t + a_2t^2 + \cdots$ \defn{normalized} if $a_1 = 1$.  The following result is due to Thoma \cite{Th}.

\begin{theorem}
The set of normalized totally positive functions is in bijection with the set of characters of the infinite symmetric group.
\end{theorem}

Via this theorem, the homomorphisms of Proposition \ref{P:hom} have the following interpretation.  Given a character $\chi: S_\infty \to \C$ we may restrict to $S_n$ to get a (reducible, fractional) character $\chi|_{S_n}$ of $S_n$.  The value $\phi(s_\lambda)$ of the corresponding homomorphism $\phi:\Sym \to \R$ is the coefficient of the irreducible character $\chi^\lambda$ in $\chi|_{S_n}$.

If a totally positive function $p(t)$ is a polynomial, then it factors into linear factors of the form $(1+\alpha t)$.  In the case that $p(t)$ is not a polynomial, the Edrei-Thoma theorem \cite{Edr,Th} gives a complete description of $p(t)$ as a meromorphic function.  In this context, Vershik and Kerov \cite{VK} give a beautiful interpretation of the poles and zeroes of $p(t)$: thinking of a character of $S_\infty$ as a point-wise limit
$$
\chi = \lim_{n \to \infty} \chi^{\lambda^{(n)}}/\chi^{\lambda^{(n)}}(1)
$$
of normalized irreducible characters of $S_n$, the poles and zeroes control the asymptotic growth of the lengths of rows and columns in $\lambda^{(n)}$.

While the following problem was one of our main motivations in \cite{TP1}, little progress has been made.
\begin{problem}
Find an interpretation of totally nonnegative matrix formal power series ($U_{\geq 0}$ for $n > 1$) in terms of representation theory of $S_\infty$, or of some other asymptotic representation theory.
\end{problem}

As Proposition \ref{P:hom} suggests, this may require expressing $\LSym$ as the character ring of a series of groups, in the way that $\Sym$ is the character ring of the series of symmetric groups.

\section{Crystals, energy, and charge}\label{sec:crystal}
In Section \ref{sec:R}, we discussed a birational action of $S_m$ on a $n \times m$ array of variables.  This action of $S_m$ turns out to be the $R$-matrix of certain products of Kirillov-Reshetikhin crystals.  We first give an informal introduction to Kashiwara's crystal graphs; focusing in particular on the relation with tableaux.

\subsection{Affine $R$-matrix and jeu de taquin}
For a combinatorial introduction to this subject, we recommend Shimozono's ``Crystals for dummies'' \cite{Sh}.  Crystal graphs were invented by Kashiwara \cite{Kas} as combinatorial models for representations of quantum groups.  For our purpose, a crystal graph (of affine type $A$) is a finite graph with edges labeled by one of $0,1,2,\ldots,n-1$.    Here is an example of an affine crystal graph:
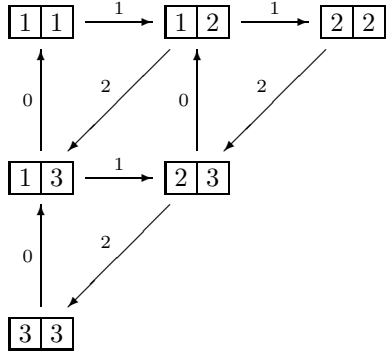
\dgARROWLENGTH=2.5em
$$
\begin{diagram}
\node{\tableau[sy]{1&1 }} \arrow{ee,t}{1} %
\node{\tableau[sy]{ 1&2}} \arrow{sw,t}{2} \arrow{ee,t}{1} %
\node{\tableau[sy]{2&2}}  \arrow{sw,t}{2} \\
\node{\tableau[sy]{1&3}} \arrow{n,t}{0} \arrow{ee,t}{1}%
\node{\tableau[sy]{2&3}} \arrow{sw,t}{2} \arrow{n,t}{0}\\
\node{\tableau[sy]{3&3}} \arrow{n,t}{0}
\end{diagram}
$$

For a positive integer $s$, the set of semistandard Young tableaux of row shape $(s)$, with entries in $1,2,\ldots,n$, form an affine crystal, which is a special case of a Kirillov-Reshetikhin crystal.

If $B_1, B_2$ are Kirillov-Reshetikhin crystals, the \defn{combinatorial $R$-matrix} is the unique isomorphism $R_{B_1,B_2}: B_1 \otimes B_2 \to B_2 \otimes B_1$ of affine crystals.  It is known that the combinatorial $R$-matrices generate an action of $S_m$ on $B_1 \otimes \cdots \otimes B_m$.

The $R$-matrix has a convenient combinatorial interpretation \cite{Sh} in terms of semistandard tableaux and the jeu de taquin algorithm \cite{EC2}. Let $b_1 \otimes b_2$ be an element of $B_1 \otimes B_2$.  Then $R_{B_1,B_2}(b_1 \otimes b_2) = c_1 \otimes c_2 \in B_2 \otimes B_1$ where $c_1, c_2$ are the unique pair of row shaped tableaux which jeu de taquin to the same tableau that $b_1$ and $b_2$ jeu de taquin to, as follows:
$$
\tableau[sY]{\bl&\bl&{1}&{1}&{3}&{3}&3 \\{1}&{2}&3} \;\; \stackrel{R_{B_1,B_2}}{\longrightarrow} \;\; \tableau[sY]{ \bl&\bl&\bl&\bl&{1}&{1}&3\\ {1}&{2}&{3}&{3}&3} \;\; \text{since both jdt to} \;\; \tableau[sY]{{1}&{1}&{1}&{3}&{3}&3 \\ 2&{3}}
$$

The action of the combinatorial $R$-matrix can be described in terms of coordinates, as follows \cite{HHIKTT}.   If $x = (x^{(1)},\ldots,x^{(n)})$ is an $n$-tuple, we denote by $\x$ the $n$-tuple $(x^{(n)},x^{(1)},\ldots,x^{(n-1)})$.  Recall that the tropicalization $\trop(f)$ of a subtraction-free rational function $f$ is obtained by replacing $+$ by $\min$, and $\times$ by $+$, and $/$ by $-$.

\begin{theorem}
Let $R_{B_1,B_2}(b_1 \otimes b_2) = c_1 \otimes c_2$ and let $x^{(r)}_1$, $x^{(r)}_2$, $y^{(r)}_1$, $y^{(r)}_2$ be the number of boxes filled with $r$-s in $b_1$, $b_2$, $c_1$, $c_2$ respectively, for $r = 1, \ldots, n$.  Then
$$
\trop(s(x_1, \x_2)) = (y_1,\bar y_2).
$$
\end{theorem}

\begin{example}
In the example above, we have $x_1 = (1,1,1)$, $\x_2 = (3,2,0)$, $y_1 = (1,1,3)$ and $\bar y_2 = (1,2,0)$.  To calculate $y_1^{(2)}$, one takes
$$
y_1^{(2)} = \trop(\x_2^{(3)}\frac{\kappa_3}{\kappa_2}) = \x_2^{3} +\min(2,4,5) - \min(2,1,3) = 0+2-1 = 1.
$$
\end{example}

\subsection{The energy function and charge}
The energy function $\D_B: B \to \Z$ of an affine crystal plays an important role from two different points of view: (1) from the representation theory point of view the energy encodes information corresponding to the affine weight of elements of the crystal, and (2) in the connection between affine crystals and vertex models, elements of a tensor product of crystals are thought of as paths, and the energy function has the physical interpretation as the energy of a path.  See \cite{KKMMNN}.

Instead of giving the usual axiomatic or recursive definitions of the energy function, we will instead describe it in relation to another well known statistic on words: \defn{cocharge}.   Suppose $u$ is a word in the letters $1,2,\ldots,n$.  Then the weight $\wt(u)$ of $u$ is the composition $\alpha=(\alpha_1,\ldots,\alpha_n)$ where $\alpha_i$ is the number of $i$'s in $u$.   Suppose the weight of $u$ is a partition, namely, $\alpha_1 \geq \alpha_2 \geq \cdots \geq \alpha_n$.  The cocharge $cc(u)$ is obtained recursively as follows.  Underline the rightmost $1$ in $u$.  Given an underlined $i-1$, underline the rightmost $i$ to its left, if it exists; otherwise, underline the rightmost $i$ in $u$, if it exists; and otherwise stop.  To each underlined letter we assign an index: the $1$ is assigned index $0$, and the index of the underlined $i$ is equal to the index of the underlined $i-1$ if the $i$ occurs to the left of the $i$, otherwise the index of the $i$ is one more than the index of the $i-1$.  Finally, we define $cc(u) = \sum \text{indices} \; + cc(u')$ where the word $u'$ is obtained from $u$ by erasing all the underlined letters.  We define the cocharge of the empty word to be 0.

\begin{example}\label{ex:cc}
The cocharge  of $u = 3222311111233$ is calculated as follows (indices are indicated as subscripts):
\begin{align*}
&\underline{3}_0 22\underline{2}_0 31111\underline{1}_0 233 \\
& 2\underline{2}_03111\underline{1}_0 23\underline{3}_1\\
&\underline{2}_0311\underline{1}_02\underline{3}_1\\
&\underline{3}_11\underline{1}_0\underline{2}_1\\
&\underline{1}_0
\end{align*}
Thus $cc(u) = (0+0+0)+(0+0+1) + (0+0+1) + (0+1+1) + 0 = 4$.
\end{example}

Let $T$ be a semistandard tableaux with partition weight $\mu_1,\mu_2,\ldots,\mu_m$.  Let $a^{(j)}_{i}(T)$ be the number of $i$'s in the $j$-th row of $T$.  We define a crystal element $b(T) = b_1 \otimes b_2 \otimes \cdots \otimes b_m \in B_{\mu_1} \otimes B_{\mu_2} \otimes \cdots \otimes B_{\mu_m}$ by setting $b_i$ to be the one-row tableau with $a^{(j)}_{i}$ $j$'s.  These $b(T)$ are exactly the highest weight vectors of the affine crystal.  Recall that the reading word $r(T)$ of a tableau is obtained by reading the rows from left to right, starting from the bottom row.  The following result is due to Nakayashiki and Yamada \cite{NY}.

\begin{theorem}
Let $T$ be a semistandard tableaux with partition weight.  Then the cocharge $cc(r(T))$ of the reading word is equal to the energy $\D_B(b(T))$ of the corresponding crystal element.
\end{theorem}

\begin{example}\label{ex:bt}
Let
$$
T = \tableau[sY]{1&1&1&1&1&2&3&3\\2&2&2&3\\3}
$$
Then $r(T)$ is the word in Example \ref{ex:cc}.  Thus the energy of
$$
b(T) =\tableau[sY]{1&1&1&1} \otimes \tableau[sY]{1&2&2&2} \otimes \tableau[sY]{1&1&2&3}
$$
is $4$.
\end{example}

It is known that the energy function of a tensor product $B_{\mu_1} \otimes B_{\mu_2} \otimes \cdots \otimes B_{\mu_m}$ of Kirillov-Reshetikhin crystals commutes with the affine $R$-matrix.  As explained in Section \ref{sec:R}, the ring of loop symmetric functions is exactly the ring of polynomial invariants of the birational $R$-matrix.  Nevertheless, the following result is somewhat surprising:

\begin{theorem}\label{T:energy}
Set $x_i^{(j)} = a_i^{(j+1-i)}$.  The energy function $\D_B$ of the affine crystal $B_{\mu_1} \otimes B_{\mu_2} \otimes \cdots \otimes B_{\mu_m}$ is equal to the tropicalization of the loop Schur function $s^{(0)}_{(m-1)\delta_{n-1}}(x_1,x_2,\ldots,x_m)$.
\end{theorem}

There seems to be no {\it a priori} reason that the energy function $\D_B$ should be the tropicalization of a {\it polynomial}, much less that of a loop Schur function.

\begin{example}
For Example \ref{ex:bt} we would have $x_1=(5,0,0)$, $x_2=(0,1,3)$, and $x_3=(1,1,2)$.  According to Theorem \ref{T:energy}, the energy is
\begin{align*}
\min(&x_1^{(3)}+x_1^{(2)}+x_1^{(1)}+x_1^{(3)}+x_2^{(1)}+x_2^{(3)}, \\ &x_1^{(3)}+x_1^{(2)}+x_2^{(1)}+x_3^{(3)}+x_3^{(1)}+x_3^{(3)}, x_1^{(3)}+x_2^{(3)}+x_2^{(1)}+x_2^{(3)}+x_2^{(1)}+x_3^{(3)},\ldots)
\end{align*}
corresponding to the tableaux
$$
\tableau[sY]{1&1&1&1\\2&2} \qquad \tableau[sY]{1&1&2&3\\3&3} \qquad \tableau[sY]{1&2&2&2\\2&3} \qquad \cdots
$$
The minimum is achieved (uniquely) on the tableau
$$
\tableau[sY]{1&1&2&3\\2&3}
$$
which gives $x_1^{(3)}+x_1^{(2)}+x_2^{(1)}+x_3^{(3)}+x_2^{(1)}+x_3^{(3)} = 0+0+0+2+0+2=4$, agreeing with the previous calculation.
\end{example}

\section{Box-ball systems}\label{sec:bb}
Loop symmetric functions can also be thought of as integrals of motion of certain discrete dynamical systems, called box-ball systems.    Here the birational action of $S_n$ is the time evolution of the system.  The connection between the $R$-matrix and the box-ball system was established in \cite{HHIKTT}.

\subsection{Takahashi-Satusuma box-ball system \cite{TS}}
We begin with an (infinite) configuration of boxes arranged on a line.  Each box can either contain one ball or no balls.  We shall always assume that there are finitely many balls.  Time evolution proceeds as follows.  Take the leftmost ball and place it in the leftmost empty box to the right of that ball.  Repeatedly do this for the leftmost ball that has not yet been moved, until all balls are moved.  This completes the time evolution.

\begin{figure}
\unitlength 15pt
\begin{picture}(20,3)
\put(1,0){\tama}
\put(2,0){\tama}
\put(3,0){\tama}
\put(8,0){\tama}
\multiput(0,0)(1,0){20}{\hako}
\end{picture}
\begin{picture}(20,3)
\put(4,0){\tama}
\put(5,0){\tama}
\put(6,0){\tama}
\put(9,0){\tama}
\multiput(0,0)(1,0){20}{\hako}
\end{picture}
\begin{picture}(20,3)
\put(7,0){\tama}
\put(8,0){\tama}
\put(10,0){\tama}
\put(11,0){\tama}
\multiput(0,0)(1,0){20}{\hako}
\end{picture}
\begin{picture}(20,3)
\put(9,0){\tama}
\put(12,0){\tama}
\put(13,0){\tama}
\put(14,0){\tama}
\multiput(0,0)(1,0){20}{\hako}
\end{picture}
\begin{picture}(20,3)
\put(10,0){\tama}
\put(15,0){\tama}
\put(16,0){\tama}
\put(17,0){\tama}
\multiput(0,0)(1,0){20}{\hako}
\end{picture}
\caption{Time evolution in the box-ball system.}
\label{fig:bb}
\end{figure}
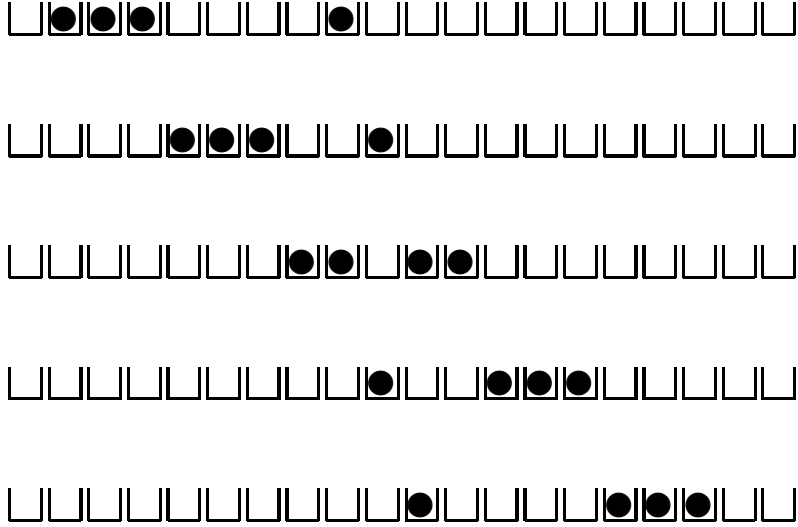

Figure \ref{fig:bb} shows such a configuration, and its time evolution.

While the definition is simple, the box-ball system has the remarkable property that it exhibits solitonic behavior.  An isolated consecutive string of $k$ balls is a \defn{soliton} which travels with speed $k$.  Figure \ref{fig:bb} shows the collision of a speed 3 solition with a speed 1 soliton.  Note how the faster soliton overtakes the slower one, but both solitons are preserved.

\subsection{Carrier description, and the birational $R$-matrix}
The time evolution of the box-ball system can also be described in terms of a carrier.  Imagine that a carrier (with infinite capacity) begins from the left going from box to box.  At each box, the carrier:
\begin{enumerate}
\item
Picks up a ball, if the box has a ball in it.
\item
Drops a ball, if the box is empty and the carrier is carrying a ball.
\end{enumerate}

It is easy to see that this carrier description gives the same outcome as the time evolution described previously.  Let us rewrite the interaction of the carrier with a site in terms of coordinates.  A site with $a$ extra capacity and $b$ balls will be denoted $(a,b)$; thus the empty box is denoted $(1,0)$, and the box with a ball in it is denoted $(0,1)$.  The carrier will be denoted $(\infty,d)$ if it has $d$ balls in it.  So when the carrier begins from the left, it is in the state $(\infty,0)$.  Then the interaction of the carrier and the box can be written as
$$
\{(a,b), (\infty,d)\} \longmapsto \{(a-\min(a,d)+b, \min(a,d)), (\infty,d-\min(a,d)+b)\}.
$$
The reader is invited to verify that this transformation is (after some relabeling) the tropicalization of the birational $S_2$-action for $n = 2$:
\begin{align*}
s(x^{(1)},x^{(2)}) &= \left(y^{(2)}\frac{x^{(1)}+y^{(1)}}{x^{(2)}+y^{(2)}}, y^{(1)}\frac{x^{(2)}+y^{(2)}}{x^{(1)}+y^{(1)}}\right)\\
s(y^{(1)},y^{(2)}) &= \left(x^{(2)}\frac{x^{(1)}+y^{(1)}}{x^{(2)}+y^{(2)}}, x^{(1)}\frac{x^{(2)}+y^{(2)}}{x^{(1)}+y^{(1)}}\right)
\end{align*}

It follows that Theorem \ref{T:fundlsym} can be interpreted as a statement concerning the integrals of motion of the box-ball system.

\begin{proposition}
The polynomial integrals of motion of the birational analogue of the box-ball system is exactly $\LSym$.
\end{proposition}

If we imagine that the system has $m$ boxes, then the invariants ``total number of balls'' and ``total amount of free space'' are equal to the tropicalizations of $e_m^{(1)}$ and $e_m^{(2)}$.  Because of the curious way the box-ball coordinates are shifted with respect to the $\LSym$-variables, the invariant $e_1^{(1)}$ can be given the following interpretation.  Its tropicalization is equal to 1 if every odd site is empty and every even site is occupied, and equal to 0 otherwise.  Each time evolution, the notions of odd and even swap (pushing the carrier through shifts the parity).  It would be interesting to give similar interpretations for a general loop (skew) Schur function.

See \cite{TTS} for a description of invariants of the box-ball system in terms of the Robinson-Schensted correspondence.

\begin{remark}
There is a vast literature generalizing box-ball systems to situations with many different kinds of balls (corresponding to our $n > 2$), or to situations where boxes have higher capacity.  See \cite{HHIKTT} for some directions.
\end{remark}
\begin{remark}\label{rem:BBBS}
Together with Pylyavskyy and Sakamoto, we give in \cite{LPS} a generalization of the box-ball system which arises from considering mixed directional horizontal wire on a cylinder, as in Section \ref{sec:net}.  This new solitonic system involves, boxes, balls, and baskets.  (Baskets can contain balls, but they are allowed to move as well.)
\end{remark}

\bibliographystyle{amsalpha}

\end{document}

%% file: wire181shift.pstex_t
\begin{picture}(0,0)%
\includegraphics{wire181shift.pstex}%
\end{picture}%
\setlength{\unitlength}{1579sp}%
\begingroup\makeatletter\ifx\SetFigFont\undefined%
\gdef\SetFigFont#1#2#3#4#5{%
  \reset@font\fontsize{#1}{#2pt}%
  \fontfamily{#3}\fontseries{#4}\fontshape{#5}%
  \selectfont}%
\fi\endgroup%
\begin{picture}(10760,2444)(3729,-4583)
\put(4426,-2611){\makebox(0,0)[lb]{\smash{{\SetFigFont{10}{12.0}{\rmdefault}{\mddefault}{\updefault}{\color[rgb]{0,0,0}$x^{(1)}$}%
}}}}
\put(4426,-3811){\makebox(0,0)[lb]{\smash{{\SetFigFont{10}{12.0}{\rmdefault}{\mddefault}{\updefault}{\color[rgb]{0,0,0}$x^{(2)}$}%
}}}}
\put(10651,-2611){\makebox(0,0)[lb]{\smash{{\SetFigFont{10}{12.0}{\rmdefault}{\mddefault}{\updefault}{\color[rgb]{0,0,0}$x'^{(1)}$}%
}}}}
\put(10651,-3811){\makebox(0,0)[lb]{\smash{{\SetFigFont{10}{12.0}{\rmdefault}{\mddefault}{\updefault}{\color[rgb]{0,0,0}$x'^{(2)}$}%
}}}}
\put(6451,-2611){\makebox(0,0)[lb]{\smash{{\SetFigFont{10}{12.0}{\rmdefault}{\mddefault}{\updefault}{\color[rgb]{0,0,0}$y^{(2)}$}%
}}}}
\put(6451,-3811){\makebox(0,0)[lb]{\smash{{\SetFigFont{10}{12.0}{\rmdefault}{\mddefault}{\updefault}{\color[rgb]{0,0,0}$y^{(1)}$}%
}}}}
\put(12751,-3811){\makebox(0,0)[lb]{\smash{{\SetFigFont{10}{12.0}{\rmdefault}{\mddefault}{\updefault}{\color[rgb]{0,0,0}$y'^{(1)}$}%
}}}}
\put(12751,-2611){\makebox(0,0)[lb]{\smash{{\SetFigFont{10}{12.0}{\rmdefault}{\mddefault}{\updefault}{\color[rgb]{0,0,0}$y'^{(2)}$}%
}}}}
\end{picture}%

%% file: wire8shrink.pstex_t
\begin{picture}(0,0)%
\includegraphics{wire8shrink.pstex}%
\end{picture}%
\setlength{\unitlength}{2368sp}%
\begingroup\makeatletter\ifx\SetFigFont\undefined%
\gdef\SetFigFont#1#2#3#4#5{%
  \reset@font\fontsize{#1}{#2pt}%
  \fontfamily{#3}\fontseries{#4}\fontshape{#5}%
  \selectfont}%
\fi\endgroup%
\begin{picture}(6917,2622)(4370,-4897)
\put(10344,-3831){\makebox(0,0)[lb]{\smash{{\SetFigFont{9}{10.8}{\rmdefault}{\mddefault}{\updefault}{\color[rgb]{0,0,0}$\frac{xy}{x+z}$}%
}}}}
\put(5374,-3034){\makebox(0,0)[lb]{\smash{{\SetFigFont{9}{10.8}{\rmdefault}{\mddefault}{\updefault}{\color[rgb]{0,0,0}$x$}%
}}}}
\put(5988,-4077){\makebox(0,0)[lb]{\smash{{\SetFigFont{9}{10.8}{\rmdefault}{\mddefault}{\updefault}{\color[rgb]{0,0,0}$y$}%
}}}}
\put(6233,-3034){\makebox(0,0)[lb]{\smash{{\SetFigFont{9}{10.8}{\rmdefault}{\mddefault}{\updefault}{\color[rgb]{0,0,0}$z$}%
}}}}
\put(8933,-3831){\makebox(0,0)[lb]{\smash{{\SetFigFont{9}{10.8}{\rmdefault}{\mddefault}{\updefault}{\color[rgb]{0,0,0}$\frac{yz}{x+z}$}%
}}}}
\put(9915,-3279){\makebox(0,0)[lb]{\smash{{\SetFigFont{9}{10.8}{\rmdefault}{\mddefault}{\updefault}{\color[rgb]{0,0,0}$x+z$}%
}}}}
\end{picture}%

%% file: wire11a.pstex_t
\begin{picture}(0,0)%
\includegraphics{wire11a.pstex}%
\end{picture}%
\setlength{\unitlength}{1579sp}%
\begingroup\makeatletter\ifx\SetFigFont\undefined%
\gdef\SetFigFont#1#2#3#4#5{%
  \reset@font\fontsize{#1}{#2pt}%
  \fontfamily{#3}\fontseries{#4}\fontshape{#5}%
  \selectfont}%
\fi\endgroup%
\begin{picture}(6419,1244)(4479,-6083)
\put(5026,-5161){\makebox(0,0)[lb]{\smash{{\SetFigFont{10}{12.0}{\rmdefault}{\mddefault}{\updefault}{\color[rgb]{0,0,0}$x$}%
}}}}
\put(6001,-5161){\makebox(0,0)[lb]{\smash{{\SetFigFont{10}{12.0}{\rmdefault}{\mddefault}{\updefault}{\color[rgb]{0,0,0}$-x$}%
}}}}
\end{picture}%

%% file: wire20shift.pstex_t
\begin{picture}(0,0)%
\includegraphics{wire20shift.pstex}%
\end{picture}%
\setlength{\unitlength}{1184sp}%
\begingroup\makeatletter\ifx\SetFigFont\undefined%
\gdef\SetFigFont#1#2#3#4#5{%
  \reset@font\fontsize{#1}{#2pt}%
  \fontfamily{#3}\fontseries{#4}\fontshape{#5}%
  \selectfont}%
\fi\endgroup%
\begin{picture}(26510,6034)(579,-5623)
\put(2401,-4186){\makebox(0,0)[lb]{\smash{{\SetFigFont{8}{9.6}{\rmdefault}{\mddefault}{\updefault}{\color[rgb]{0,0,0}.}%
}}}}
\put(2401,-4411){\makebox(0,0)[lb]{\smash{{\SetFigFont{8}{9.6}{\rmdefault}{\mddefault}{\updefault}{\color[rgb]{0,0,0}.}%
}}}}
\put(2401,-4636){\makebox(0,0)[lb]{\smash{{\SetFigFont{8}{9.6}{\rmdefault}{\mddefault}{\updefault}{\color[rgb]{0,0,0}.}%
}}}}
\put(17701,-4786){\makebox(0,0)[lb]{\smash{{\SetFigFont{8}{9.6}{\rmdefault}{\mddefault}{\updefault}{\color[rgb]{0,0,0}.}%
}}}}
\put(17701,-5011){\makebox(0,0)[lb]{\smash{{\SetFigFont{8}{9.6}{\rmdefault}{\mddefault}{\updefault}{\color[rgb]{0,0,0}.}%
}}}}
\put(17701,-5236){\makebox(0,0)[lb]{\smash{{\SetFigFont{8}{9.6}{\rmdefault}{\mddefault}{\updefault}{\color[rgb]{0,0,0}.}%
}}}}
\put(17476,-136){\makebox(0,0)[lb]{\smash{{\SetFigFont{8}{9.6}{\rmdefault}{\mddefault}{\updefault}{\color[rgb]{0,0,0}$-p$}%
}}}}
\put(20701,-2461){\makebox(0,0)[lb]{\smash{{\SetFigFont{8}{9.6}{\rmdefault}{\mddefault}{\updefault}{\color[rgb]{0,0,0}.}%
}}}}
\put(20926,-2461){\makebox(0,0)[lb]{\smash{{\SetFigFont{8}{9.6}{\rmdefault}{\mddefault}{\updefault}{\color[rgb]{0,0,0}.}%
}}}}
\put(21151,-2461){\makebox(0,0)[lb]{\smash{{\SetFigFont{8}{9.6}{\rmdefault}{\mddefault}{\updefault}{\color[rgb]{0,0,0}.}%
}}}}
\put(24601,-4936){\makebox(0,0)[lb]{\smash{{\SetFigFont{8}{9.6}{\rmdefault}{\mddefault}{\updefault}{\color[rgb]{0,0,0}.}%
}}}}
\put(24601,-5161){\makebox(0,0)[lb]{\smash{{\SetFigFont{8}{9.6}{\rmdefault}{\mddefault}{\updefault}{\color[rgb]{0,0,0}.}%
}}}}
\put(24601,-5386){\makebox(0,0)[lb]{\smash{{\SetFigFont{8}{9.6}{\rmdefault}{\mddefault}{\updefault}{\color[rgb]{0,0,0}.}%
}}}}
\put(24376,-1486){\makebox(0,0)[lb]{\smash{{\SetFigFont{8}{9.6}{\rmdefault}{\mddefault}{\updefault}{\color[rgb]{0,0,0}$-p$}%
}}}}
\put(24526, 14){\makebox(0,0)[lb]{\smash{{\SetFigFont{8}{9.6}{\rmdefault}{\mddefault}{\updefault}{\color[rgb]{0,0,0}$p$}%
}}}}
\put(3151,-811){\makebox(0,0)[lb]{\smash{{\SetFigFont{8}{9.6}{\rmdefault}{\mddefault}{\updefault}{\color[rgb]{0,0,0}$y^{(2)}$}%
}}}}
\put(3151,-2011){\makebox(0,0)[lb]{\smash{{\SetFigFont{8}{9.6}{\rmdefault}{\mddefault}{\updefault}{\color[rgb]{0,0,0}$y^{(3)}$}%
}}}}
\put(3151,-3211){\makebox(0,0)[lb]{\smash{{\SetFigFont{8}{9.6}{\rmdefault}{\mddefault}{\updefault}{\color[rgb]{0,0,0}$y^{(4)}$}%
}}}}
\put(7501,-4936){\makebox(0,0)[lb]{\smash{{\SetFigFont{8}{9.6}{\rmdefault}{\mddefault}{\updefault}{\color[rgb]{0,0,0}.}%
}}}}
\put(7501,-5161){\makebox(0,0)[lb]{\smash{{\SetFigFont{8}{9.6}{\rmdefault}{\mddefault}{\updefault}{\color[rgb]{0,0,0}.}%
}}}}
\put(7501,-5386){\makebox(0,0)[lb]{\smash{{\SetFigFont{8}{9.6}{\rmdefault}{\mddefault}{\updefault}{\color[rgb]{0,0,0}.}%
}}}}
\put(7426,-1561){\makebox(0,0)[lb]{\smash{{\SetFigFont{8}{9.6}{\rmdefault}{\mddefault}{\updefault}{\color[rgb]{0,0,0}$p$}%
}}}}
\put(7276, 14){\makebox(0,0)[lb]{\smash{{\SetFigFont{8}{9.6}{\rmdefault}{\mddefault}{\updefault}{\color[rgb]{0,0,0}$-p$}%
}}}}
\put(12601,-4786){\makebox(0,0)[lb]{\smash{{\SetFigFont{8}{9.6}{\rmdefault}{\mddefault}{\updefault}{\color[rgb]{0,0,0}.}%
}}}}
\put(12601,-5011){\makebox(0,0)[lb]{\smash{{\SetFigFont{8}{9.6}{\rmdefault}{\mddefault}{\updefault}{\color[rgb]{0,0,0}.}%
}}}}
\put(12601,-5236){\makebox(0,0)[lb]{\smash{{\SetFigFont{8}{9.6}{\rmdefault}{\mddefault}{\updefault}{\color[rgb]{0,0,0}.}%
}}}}
\put(12376,-136){\makebox(0,0)[lb]{\smash{{\SetFigFont{8}{9.6}{\rmdefault}{\mddefault}{\updefault}{\color[rgb]{0,0,0}$-p$}%
}}}}
\put(12526,-2311){\makebox(0,0)[lb]{\smash{{\SetFigFont{8}{9.6}{\rmdefault}{\mddefault}{\updefault}{\color[rgb]{0,0,0}$p'$}%
}}}}
\put(8251,-1561){\makebox(0,0)[lb]{\smash{{\SetFigFont{8}{9.6}{\rmdefault}{\mddefault}{\updefault}{\color[rgb]{0,0,0}$y^{(2)}$}%
}}}}
\put(8251,-2761){\makebox(0,0)[lb]{\smash{{\SetFigFont{8}{9.6}{\rmdefault}{\mddefault}{\updefault}{\color[rgb]{0,0,0}$y^{(3)}$}%
}}}}
\put(8251,-3961){\makebox(0,0)[lb]{\smash{{\SetFigFont{8}{9.6}{\rmdefault}{\mddefault}{\updefault}{\color[rgb]{0,0,0}$y^{(4)}$}%
}}}}
\put(13351,-811){\makebox(0,0)[lb]{\smash{{\SetFigFont{8}{9.6}{\rmdefault}{\mddefault}{\updefault}{\color[rgb]{0,0,0}$y'^{(2)}$}%
}}}}
\put(13351,-3811){\makebox(0,0)[lb]{\smash{{\SetFigFont{8}{9.6}{\rmdefault}{\mddefault}{\updefault}{\color[rgb]{0,0,0}$y^{(4)}$}%
}}}}
\put(18451,-811){\makebox(0,0)[lb]{\smash{{\SetFigFont{8}{9.6}{\rmdefault}{\mddefault}{\updefault}{\color[rgb]{0,0,0}$y'^{(2)}$}%
}}}}
\put(18451,-3811){\makebox(0,0)[lb]{\smash{{\SetFigFont{8}{9.6}{\rmdefault}{\mddefault}{\updefault}{\color[rgb]{0,0,0}$y^{(4)}$}%
}}}}
\put(25351,-1561){\makebox(0,0)[lb]{\smash{{\SetFigFont{8}{9.6}{\rmdefault}{\mddefault}{\updefault}{\color[rgb]{0,0,0}$y'^{(2)}$}%
}}}}
\put(25351,-2761){\makebox(0,0)[lb]{\smash{{\SetFigFont{8}{9.6}{\rmdefault}{\mddefault}{\updefault}{\color[rgb]{0,0,0}$y'^{(3)}$}%
}}}}
\put(25351,-3961){\makebox(0,0)[lb]{\smash{{\SetFigFont{8}{9.6}{\rmdefault}{\mddefault}{\updefault}{\color[rgb]{0,0,0}$y'^{(4)}$}%
}}}}
\put(18451,-2011){\makebox(0,0)[lb]{\smash{{\SetFigFont{8}{9.6}{\rmdefault}{\mddefault}{\updefault}{\color[rgb]{0,0,0}$y'^{(3)}$}%
}}}}
\put(13351,-2611){\makebox(0,0)[lb]{\smash{{\SetFigFont{8}{9.6}{\rmdefault}{\mddefault}{\updefault}{\color[rgb]{0,0,0}$y^{(3)}$}%
}}}}
\put(1051,-811){\makebox(0,0)[lb]{\smash{{\SetFigFont{8}{9.6}{\rmdefault}{\mddefault}{\updefault}{\color[rgb]{0,0,0}$x^{(1)}$}%
}}}}
\put(1051,-2011){\makebox(0,0)[lb]{\smash{{\SetFigFont{8}{9.6}{\rmdefault}{\mddefault}{\updefault}{\color[rgb]{0,0,0}$x^{(2)}$}%
}}}}
\put(1051,-3211){\makebox(0,0)[lb]{\smash{{\SetFigFont{8}{9.6}{\rmdefault}{\mddefault}{\updefault}{\color[rgb]{0,0,0}$x^{(3)}$}%
}}}}
\put(6151,-1561){\makebox(0,0)[lb]{\smash{{\SetFigFont{8}{9.6}{\rmdefault}{\mddefault}{\updefault}{\color[rgb]{0,0,0}$x^{(1)}$}%
}}}}
\put(6151,-2761){\makebox(0,0)[lb]{\smash{{\SetFigFont{8}{9.6}{\rmdefault}{\mddefault}{\updefault}{\color[rgb]{0,0,0}$x^{(2)}$}%
}}}}
\put(6151,-3961){\makebox(0,0)[lb]{\smash{{\SetFigFont{8}{9.6}{\rmdefault}{\mddefault}{\updefault}{\color[rgb]{0,0,0}$x^{(3)}$}%
}}}}
\put(11176,-811){\makebox(0,0)[lb]{\smash{{\SetFigFont{8}{9.6}{\rmdefault}{\mddefault}{\updefault}{\color[rgb]{0,0,0}$x'^{(1)}$}%
}}}}
\put(11176,-2611){\makebox(0,0)[lb]{\smash{{\SetFigFont{8}{9.6}{\rmdefault}{\mddefault}{\updefault}{\color[rgb]{0,0,0}$x^{(2)}$}%
}}}}
\put(11176,-3811){\makebox(0,0)[lb]{\smash{{\SetFigFont{8}{9.6}{\rmdefault}{\mddefault}{\updefault}{\color[rgb]{0,0,0}$x^{(3)}$}%
}}}}
\put(16276,-811){\makebox(0,0)[lb]{\smash{{\SetFigFont{8}{9.6}{\rmdefault}{\mddefault}{\updefault}{\color[rgb]{0,0,0}$x'^{(1)}$}%
}}}}
\put(16276,-2011){\makebox(0,0)[lb]{\smash{{\SetFigFont{8}{9.6}{\rmdefault}{\mddefault}{\updefault}{\color[rgb]{0,0,0}$x'^{(2)}$}%
}}}}
\put(16276,-3811){\makebox(0,0)[lb]{\smash{{\SetFigFont{8}{9.6}{\rmdefault}{\mddefault}{\updefault}{\color[rgb]{0,0,0}$x^{(3)}$}%
}}}}
\put(23176,-1561){\makebox(0,0)[lb]{\smash{{\SetFigFont{8}{9.6}{\rmdefault}{\mddefault}{\updefault}{\color[rgb]{0,0,0}$x'^{(1)}$}%
}}}}
\put(23176,-2761){\makebox(0,0)[lb]{\smash{{\SetFigFont{8}{9.6}{\rmdefault}{\mddefault}{\updefault}{\color[rgb]{0,0,0}$x'^{(2)}$}%
}}}}
\put(23176,-3961){\makebox(0,0)[lb]{\smash{{\SetFigFont{8}{9.6}{\rmdefault}{\mddefault}{\updefault}{\color[rgb]{0,0,0}$x'^{(3)}$}%
}}}}
\end{picture}%

%% file: wire2shrink.pstex_t
\begin{picture}(0,0)%
\includegraphics{wire2shrink.pstex}%
\end{picture}%
\setlength{\unitlength}{1973sp}%
\begingroup\makeatletter\ifx\SetFigFont\undefined%
\gdef\SetFigFont#1#2#3#4#5{%
  \reset@font\fontsize{#1}{#2pt}%
  \fontfamily{#3}\fontseries{#4}\fontshape{#5}%
  \selectfont}%
\fi\endgroup%
\begin{picture}(8793,3348)(1236,-4677)
\put(3384,-1997){\makebox(0,0)[lb]{\smash{{\SetFigFont{8}{9.6}{\rmdefault}{\mddefault}{\updefault}{\color[rgb]{0,0,0}$f_2$}%
}}}}
\put(1333,-1997){\makebox(0,0)[lb]{\smash{{\SetFigFont{8}{9.6}{\rmdefault}{\mddefault}{\updefault}{\color[rgb]{0,0,0}$e_1$}%
}}}}
\put(1333,-3983){\makebox(0,0)[lb]{\smash{{\SetFigFont{8}{9.6}{\rmdefault}{\mddefault}{\updefault}{\color[rgb]{0,0,0}$e_2$}%
}}}}
\put(3384,-3983){\makebox(0,0)[lb]{\smash{{\SetFigFont{8}{9.6}{\rmdefault}{\mddefault}{\updefault}{\color[rgb]{0,0,0}$f_1$}%
}}}}
\put(2524,-2990){\makebox(0,0)[lb]{\smash{{\SetFigFont{8}{9.6}{\rmdefault}{\mddefault}{\updefault}{\color[rgb]{0,0,0}$x$}%
}}}}
\put(6427,-2329){\makebox(0,0)[lb]{\smash{{\SetFigFont{8}{9.6}{\rmdefault}{\mddefault}{\updefault}{\color[rgb]{0,0,0}$x$}%
}}}}
\put(7752,-2329){\makebox(0,0)[lb]{\smash{{\SetFigFont{8}{9.6}{\rmdefault}{\mddefault}{\updefault}{\color[rgb]{0,0,0}$y$}%
}}}}
\put(8677,-3387){\makebox(0,0)[lb]{\smash{{\SetFigFont{8}{9.6}{\rmdefault}{\mddefault}{\updefault}{\color[rgb]{0,0,0}$z$}%
}}}}
\end{picture}%